\documentclass[12pt]{article}

\title{Some geometric aspects of Puiseux surfaces}
\author{Jos\'e M. Tornero}
\date{31/10/2002}

\newcommand{\CC}{{\mathbf C}}
\newcommand{\ZZ}{{\mathbf Z}}
\newcommand{\NN}{{\mathbf N}}

\newcommand{\cS}{\mathcal S}
\newcommand{\ol}{\overline}

\begin{document}

\maketitle

\begin{abstract}
The following problem is treated: Characterizing the tangent
cone and the equimultiple locus of a Puiseux surface (that
is, an algebroid embedded surface admitting an equation whose roots are Puiseux
power series), using a set of exponents appearing in a root of an equation.
The aim is knowing to which extent the well--known results for the
quasi--ordinary case can be extended to this much wider family.
\end{abstract}

\vspace{.3cm}



\vspace{.3cm}

\section{Introduction. Privileged exponents}

Let us consider $\CC$ (or any other algebraically closed field of
characteristic $0$), let $X$ and $Y$ be formally independent variables over
$\CC$ and let $m \in \NN$ be a positive integer, fixed in what follows. We
have the following rings and fields:
$$
K = \CC ((X,Y)) \subset L = \CC \left( \left( X^{1/m},Y^{1/m}
\right) \right)
$$
$$
R = \CC [[X,Y]] \subset S = \CC \left[ \left[ X^{1/m},Y^{1/m}
\right] \right].
$$
The elements of $S$ will be called Puiseux power series.

The field extension is trivially normal and finite, having as Galois group
$G \simeq C_m \times C_m$. We will denote the elements of $G$ by
\begin{eqnarray*}
(a,b): L & \longrightarrow & L \\
X^{1/m} & \longmapsto & \delta^a X^{1/m} \\
Y^{1/m} & \longmapsto & \delta^b Y^{1/m}
\end{eqnarray*} 
where $\delta$ (fixed from now on) is a primitive $m$--th root of the unity.

\noindent {\bf Definition.--}
Let $\cS=Spec(\CC[[X,Y,Z]]/(F)) = Spec (R_\cS)$ be an embedded algebroid 
surface. Then $\cS$ is called a Puiseux surface if the roots of $F$ are
Puiseux power series.

If $\cS$ is an irreducible Puiseux surface and $F$ is an equation of 
$\cS$ which is the minimal polynomial of $\zeta \in S$, we will say that
$\cS$ is defined by $\zeta$.

\noindent {\bf Definition.--}
Let $\zeta \in S$, written as
$$
\zeta = \sum_{(i,j) \in \Delta (\zeta) \subset \NN^2} c_{ij}
X^{i/m} Y^{j/m}, 0 \neq c_{ij} \in \CC.
$$

The set $\Delta (\zeta)$ ($\Delta$, if no confussion is posible) is called the 
set of exponents of $\zeta$. Note that we consider the exponents $(i,j)$ rather
than $(i/m,j/m)$.

A particular case of Puiseux power series (and the surfaces they determine) has
been studied quite thoroughly: that of the quasi--ordinary series. A power
series $\zeta \in S$ with minimal polynomial $F(Z)$ is called quasi--ordinary
if
$$
D(F) = X^a Y^b u(X,Y), \mbox{ with } u(0,0) \neq 0,
$$
where $D(F)$ is the discriminant of $F$ (with respect to $Z$, of course). These
series were already used by Jung (\cite{Jung}) in his work on the resolution of
surface singularities. Later on, the work of Lipman (\cite{Lipman}) and others
(see, for instance \cite{Herrera}, \cite{Gau}) has led to a quite well
understanding of the geometry and the topology of these surfaces and their
resolution.

In particular, it is known that the quasi--ordinary character is preserved by
blowing--up equimultiple smooth subvarieties.

\noindent {\bf Definition.--}
Given a Puiseux power series $\zeta = \sum c_{ij} X^{i/m}Y^{j/m}$, we will say
that a set $E=\{(i_1,j_1),...,(i_t,j_t)\} \subset \Delta$ is a set of
privileged exponents if
$$
K(\zeta) = K \left( X^{i_1/m}Y^{j_1/m},...,X^{i_t/m}Y^{j_t/m} \right).
$$

\noindent {\bf Example.--}
Consider the Puiseux power series
$$
\zeta = X^{13/9}Y^{16/9} - 2 X^{22/9}Y^{7/9} + 7 X^{12/9} Y^{11/9} -5X^{7/9}
Y^{4/9}.
$$

Then we have the following minimal sets of privileged exponents:
$$
E_1= \{ (7,4), \, (12,11) \}, \; E_2= \{ (12,11), \, (13, 16) \}, \; E_3 = \{ (12,11), \,
(22,7) \}.
$$

In fact, these are all the minimal sets of privileged exponents, as
$$
X^{12/9}Y^{11/9} \notin K \left( X^{13/9}Y^{16/9},X^{22/9}Y^{7/9},X^{7/9}Y^{4/9} 
\right).
$$

It is possible to give a constructive (in an ample sense, taking into account
that our objects are power series) process which produces a set of privileged 
exponents from a given Puiseux power series (\cite{Kummer}). The process relies
on a chosen monomial order. If a graded order is taken, the exponents for 
a quasi--ordinary power series are precisely the so--called characteristic
exponents (\cite{Lipman}) of the power series.

Let $\zeta$ be a given quasi--ordinary power series $\zeta$ with minimal 
polynomial $F$, defining a Puiseux surface $\cS = Spec (R_\cS)$. The
characteristic exponents of $\zeta$ have the following properties (see the 
above references), up to a technical process called normalization:

\begin{enumerate}
\item[(a)] They are totally ordered by the product ordering
$$
(i,j) < \left(i',j' \right) \; \Longleftrightarrow \; i<i', \; j<j'.
$$
\item[(b)] Every exponent $(i,j) \in \Delta$ can be written (mod $\ZZ m \times
\ZZ m$) as an integral combination of the characteristic exponents.
\item[(c)] Their monomials generate the same field extension (of $K$) than 
$\zeta$.
\item[(d)] They determine the degree of $F$ (that is, $[K(\zeta):K]$).
\item[(e)] They determine the set of equimultiple subvarieties of $\cS$.
\item[(f)] They determine the tangent cone of $\cS$.
\item[(g)] They determine ({\em and} are determined by) the topological type
of $\cS$ (in the sense of \cite{SE}).
\item[(h)] They determine the characteristic exponents of any permissible
blowing--up of $\cS$.
\end{enumerate}

Our abstract aim is studying to which extent the properties of quasi--ordinary
singularities can be extended to the general Puiseux singularities. In 
particular we wanted to know if the tangent cone and the equimultiple locus can
be obtained from a finite set of data. 

\noindent {\bf Remark.--}
When working with a given $\zeta \in S$ and its minimal polynomial $F(Z)$, we
can always consider that $\Delta \cap (\ZZ m)^2 = \emptyset$. This can 
be achieved by performing the Tchirnahusen transformation on $F$, that is, if
$$
F(Z) = Z^n + a_{n-1}(X,Y) Z^{n-1}+...+a_0 (X,Y),
$$
then the change of variables
$$
Z \; \longmapsto Z - \frac{1}{n} a_{n-1} (X,Y)
$$
will take $F$ to the minimal polynomial of the power series resulting from
erasing all the monomials in $\zeta$ with exponents in $(\ZZ m)^2$. So, from
now on, we will assume that all the power series we are working with have no
exponents in $(\ZZ m)^2$.

Let $\zeta \in S$ and $\{ (i_1,j_1),...,(i_t,j_t)\}$ be a set of privileged
exponents of $\zeta$.
We have already proved in \cite{Kummer} that this set verify the properties
(b), (c) (well, this is obvious) and (d). In particular, 
with the notations above
$$
\left[ K (\zeta): K \right] = \frac{m^2}{d(\zeta)},
$$
where $d(\zeta)$ is the $\gcd$ of the $2 \times 2$ minors of the matrix
$$
\left( \begin{array}{ccccc}
m & 0 & i_1 & ... & i_t \\
0 & m & j_1 & ... & j_t \\
\end{array} \right)
$$
which does not depend on the set of privileged exponents considered.

In this paper we prove that it is not possible, with full generality, to obtain 
(e) and (f) from a subset of $\Delta$. As the tangent cone and the equimultiple
subvarieties are 
essential for a complete understanding of the resolution process, we conclude
that a more complicated set of data is necessary for the general Puiseux case. 

\section{The tangent cone of Puiseux surfaces}

For our purposes, we can take the tangent cone of an embedded algebroid 
surface to be the affine variety defined by the initial form of an equation of
the surface. In fact, the actual definition is the spectrum of the graded ring
gr$_{(X,Y,Z)}(R_\cS)$, but it is easily shown that both versions
coincide (\cite{Lipman}).

\noindent {\bf Remark.--}
We can restrict ourselves to irreducible Puiseux surfaces, as the tangent cone
of a reducible surface is the union of the tangent cones of its irreducible 
components.

Let therefore $\cS$ be an irreducible Puiseux surface, defined by $\zeta$. As 
an equation of $\cS$ is
$$
F(Z) = \prod_{i=1}^n \left( Z - \zeta_i \right),
$$
where $\{ \zeta = \zeta_1,...,\zeta_n \} = \{ (a,b)(\zeta) \; | \; (a,b) \in
G\}$, it is plain from \cite{Kummer} that the multiplicity of $\cS$ is 
determined by a set of privileged exponents (for instance, choosing any graded
ordering in $\NN^2$ so that the first exponent determines the order of 
$\zeta$). In fact,
$$
mult (\cS) = \min \{ n,n \nu (\zeta) \} = \min \left\{ \frac{m^2}{d(\zeta)},
\frac{m^2 \nu (\zeta)}{d(\zeta)} \right\},
$$
where $\nu$ is the usual order in $S$. Here we need to make a crucial 
distinction for the sequel.

\noindent {\bf Definition.--}
A Puiseux power series $\zeta \in S$ will be called transversal if $\nu (\zeta)
\geq 1$. Otherwise it will be called non--transversal.

In geometrical terms, transversality (respectively, non--transversality) 
corresponds to the situation in which the generic fiber of the projection
$$
(X,Y,Z) \longmapsto (X,Y)
$$ 
has exactly $n= [K(\zeta):K]$ (respectively strictly less) points.

\subsection{The transversal case}

This case is the easiest one. Note that, if $\nu(\zeta) \geq 1$, then its
minimal polynomial
$$
F(Z) = Z^n + \sum_{l=0}^{n-2} a_l (X,Y) Z^l
$$
is a Weierstrass polynomial, that is $\nu (a_l) \geq n-l$.

\noindent {\bf Remark.--}
If $\nu(\zeta) > 1$, the tangent cone of $\cS$ is the plane $Z^n=0$.

\noindent {\bf Remark.--}
If $\nu (\zeta) = 1$, let $\ol{\zeta}$ be the initial form of $\zeta$. As
Galois conjugation leaves the exponents unchanged, we have
$$
\ol{F} = \prod_{i=1}^n \left( Z - \ol{\zeta_i} \right)
$$
with $\Delta (\zeta_i) = \Delta (\zeta)$ and $\Delta (\ol{\zeta}) = \Delta
(\ol{\zeta_i})$.

Obviously there can be repeated factors in the above decomposition of $\ol{F}$.
In particular, it becomes clear that $\ol{F}$ must be a power of the minimal
polynomial of $\ol{\zeta}$.

Using the previous expression for the degree of a minimal polynomial from a
set of privileged exponents, we can conclude that, in this case, $\ol{F}$ 
is a $d (\ol{\zeta}) / d (\zeta) $ power of an irreducible polynomial (which
is the equation of the surface defined by $\ol{\zeta}$).

\subsection{The non--transversal case}

Let us write $\nu (\zeta) = \lambda = \mu / n < 1$. In this case the situation
is less straightforward. The initial form of $F$ is now
$$
\ol{F} = \ol{\zeta_1}...\ol{\zeta_n},
$$
which is a form of order $\lambda m^2 / d (\zeta)$. Being a homogeneous 
polynomial in two variables with coefficients in $\CC$ it can be written as
a product of linear forms
$$
\ol{F} = G_1^{a_1}...G_r^{a_r}, \; \mbox{ with } G_i \neq G_j,
$$
that is, the tangent cone is a union of linear planes. It could well be
an only plane or a proper product of planes.

Let us begin with a particularly simple case:

\noindent {\bf Lemma.--}
Let $\zeta = X^{1/m} + \alpha Y^{1/m}$, $0 \neq \alpha \in \CC$. Then, the
tangent cone of the surface defined by $\zeta$ is given by $\left( X -
\alpha^m Y \right)^m$.

\noindent {\bf Proof:}
Apply Cardano formulae for the cyclotomic case.

In fact, it is easy showing that
$$
\prod_{\delta \; | \; \delta^m = 1} \left( X^{1/m} + \delta \alpha Y^{1/m}
\right) = X - \alpha^m Y.
$$

The factors on the right hand side will be called the monic conjugates of
$\zeta$.

Let us move on to the general case and let $\zeta \in S$ with
$$
\ol{\zeta} = L^{(1)}_1...L^{(1)}_{r_1}L^{(2)}_1...L^{(2)}_{r_2}...
L^{(l)}_1...L^{(l)}_{r_l},
$$
where
$$
L^{(a)}_b = X^{1/m} - \alpha^{(a)}_b Y^{1/m},
$$
and they are written in such a way that two linear forms $L^{(a)}_{b_1}$ and
$L^{(a)}_{b_2}$ share the same superindex if and only there exists $\delta$, 
an $m$--th root of the unity, with
$$
\delta \alpha^{(a)}_{b_1} = \alpha^{(a)}_{b_2}.
$$

Clearly $r_1+...+r_l = \mu$.

Moreover, when conjugating $\zeta$ by $(a,b) \in G$, $L_1^{(1)}$ must be taken 
(up to product by a constant) into one of its monic conjugates, and so does
all $L^{(1)}_j$, $j=2,...,r_1$. In fact, when we have computed all the
$m^2/d(\zeta)$ posible conjugates of $\zeta$ appearing in the decomposition of
$F$, the $m$ monic conjugates of $L_1^{(1)}$ must appear among all 
the conjugates of the $L^{(1)}_j$'s

Further, all of these monic conjugates must appear the same number of times. If
this were not so, after gathering together the remaining forms, we would get
a product
$$
\prod_{i=1}^p \left( X^{1/m} - \alpha_i Y^{1/m} \right) = s(X,Y) \in R.
$$

Now, as this product must remain invariant by any $(a,b) \in G$, the unique
factorization in $S$ implies that all the monic conjugates of every factor
must lie in the product.

Therefore, the total number of monic conjugates of the $L^{(1)}_j$'s appearing
in the constant term of $F$ is $m^2r_1/d(\zeta)$, and all of them
appear the same number of times. So, together with the previous lemma,
this proves that the decomposition of the tangent cone is
$$
\left( X - \alpha_1^m Y \right)^{mr_1/d} ...
\left( X - \alpha_l^m Y \right)^{mr_l/d}.
$$

\noindent {\bf Example.--}
If we consider
$$
\zeta = X^{3/4} + 2X^{2/4}Y^{1/4} - X^{1/4}Y^{2/4} - 2Y^{3/4} + X^{6/4},
$$
we have the following decomposition for $\ol{\zeta}$:
$$
\begin{array}{lll}
\bar{\zeta} & = & X^{3/4} + 2X^{2/4}Y^{1/4} - X^{1/4}Y^{2/4} - 2Y^{3/4} \\
& = & \left( X^{1/4} - Y^{1/4} \right) \left( X^{1/4} + Y^{1/4} \right)
\left( X^{1/4} + 2Y^{1/4} \right).
\end{array}
$$

In the previous notation, this means
$$
\begin{array}{lll}
L^{(1)}_1 = X^{1/4} - Y^{1/4}, & L^{(1)}_2 = X^{1/4} + Y^{1/4}, & r_1=2 \\
\\
L^{(2)}_1 = X^{1/4} + 2Y^{1/4}, & & r_2 = 1
\end{array}
$$
and so the decomposition of the tangent cone is
$$
\left[ (X - Y)^4 \right]^2 \left[ (X - 2^4 Y)^4 \right]^1 =
\left( X - Y \right)^8 (X-16Y)^4.
$$

\noindent {\bf Remark.--} 
In the quasi--ordinary case, Lipman showed (\cite{Lipman}) that the tangent
cone could consist in, at most, two planes. In our case, there is
no bound for the amount of planes appearing, as we will show now.

Obviously, for $m$ fixed, at most $m-1$ planes can occur to join in the
tangent cone, but there can be in fact $m-1$ of them. A simple example is
obtained by taking
$$
\zeta = \ol{\zeta} = \prod_{k=1}^{\mu} \left( X^{1/m} + m^k Y^{1/m} \right),
$$
whose tangent cone is the union of $\mu$ differents planes. By varying
$\mu$ from $1$ to $m-1$ we obtain all the possible structures.

We summarize the results on this section.

\noindent {\bf Proposition.--}
Let $\cS$ be a Puiseux surface defined by $\zeta \in S$ and $F(Z)$ the minimal
polynomial of $\zeta$; $n$ and $d(\zeta)$ defined as above. Then

\begin{enumerate}
\item[(a)] If $\nu(\zeta)>1$ the equation of the tangent cone is the $Z^n=0$.
\item[(b)] If $\nu(\zeta)=1$ the equation of the tangent cone is the
$d(\ol{\zeta})/ d(\zeta)$--th power of an irreducible polynomial (the minimal 
polynomial of $\ol{\zeta}$).
\item[(c)] If $\nu(\zeta)<1$ the equation of the tangent cone is a product of 
linear forms.
\end{enumerate}

\noindent {\bf Remark.--}
Note that, while $\nu(\zeta)$ can be known from a set of privileged 
exponents, the number of different linear forms occuring in
case (c) cannot be found out from such a set. A simple example is given by
$$
\zeta_1 = X^{2/4} + (1+i)X^{1/4}Y^{1/4} + iY^{2/4}, \; \;
\zeta_2 = X^{2/4} + 2X^{1/4}Y^{1/4} + 2Y^{2/4}
$$
which have the same $\Delta$ (hence the same set of privileged exponents for
any choice of a total ordering in $\NN^2$). However, the first power series 
defines a surface whose tangent cone is given by $(X-Y)^4$, while the tangent
cone of the surface defined by the second is $(X-Y)^2(X-4Y)^2$.

\section{The equimultiple locus of Puiseux surfaces}

The equimultiple locus of an algebroid surface is the set of primes $P \subset R[[Z]]$
such that  $F \in P^{(\nu)}$, where $\nu$ is the usual order of $F$.  We will be 
interested only on primes of height $1$ as the maximal ideal of $R[[Z]]$ lies always
in the equimultiple locus.

\noindent {\bf Remark.--}
Analogously as in the previous section, we will reduce ourselves to the
irreducible case, as a given prime is equimultiple in an algebroid surface
if and only if it is equimultiple in all of its irreducible components.

\subsection{The transversal case}

Again this is the easiest part. Firstly note that, due to the fact that the minimal 
polynomial has the form
$$
F(Z) = Z^n + a_{n-2} (X,Y) Z^{n-2} + ... + a_0 (X,Y),
$$
with $\nu(F) = n$, a prime ideal belonging to the equimultiple locus can be
taken to be $(Z,G(X,Y))$,with $G(X,Y) \in R$, $\nu (G) > 0$.

Let us begin by considering the simplest
posible question: is $(X,Z)$ in the equimultiple locus? If $\cS$ is defined by
$\zeta$ and the corresponding equation is
$$
F(Z) = \prod_{i=1}^n \left( Z - \zeta_i \right) = Z^n + \sum_{l=0}^{n-2}
a_l (X,Y) Z^l,
$$
with $\nu (a_l) \geq n-l$, it becomes obvious that $(X,Z)$ lies in the singular
locus if and only if $X^{n-l} | a_l$ for all $l=0,...,n-2$.

But every monomial of $a_l$ is of the form
$$
\alpha X^{(i_1+...+i_l)/m}Y^{(j_1+...+j_l)/m},
$$
for some choice $\{ (i_1,j_1),...,(i_l,j_l) \} \subset \Delta$. In
particular, if we choose $(i_0,j_0) \in \Delta$ and take $i_1=...=i_l=i_0$ and
$j_1=...=j_l=j_0$ we will get the monomial
$$
\alpha X^{li_0/m}Y^{lj_0/m}.
$$

Hence, it must be $li_0/m > l$, that is, $i_0>m$ for all $(i_0,j_0) \in
\Delta$. This condition clearly suffices, but it is also necessary. For, if we 
had an exponent $(i_0,j_0) \in \Delta$ with $i_0<m$, the $n$--th power of this
exponent will appear as an exponent in $a_0(X,Y)$ (at least, if we take 
$(i_0,j_0)$ with minimal degree in $X$ for avoiding cancellations) and its
degree in $X$ will be strictly lower than $n$, hence $(X,Z)$ cannot be in
the equimultiple locus.

Similarly, $(Y,Z)$ lies in the equimultiple locus if and only if $Y | \zeta$. Note
that this condition can be observed in a set of privileged exponents if
we choose, for instance, the inverse lexicographic ordering as total ordering
in $\NN^2$ (see \cite{Kummer} for how this implies $Y|\zeta$).

The general case follows the same lines: a prime $(Z,c(X,Y))$ lies in the
equimultiple locus if and only if $c(X,Y)^{n-l}|a_l$, for all $l=0,...,n-2$. So,
it obvioulsy suffices that $c(X,Y) | \zeta$.

Let us see that this is also necessary. First $c(X,Y) | \zeta$ is equivalent
to $c(X,Y) | \zeta_i$ for all $i=1,...,n$, as $c(X,Y) \in R$. So if $c(X,Y)$
lies in the equimultiple locus but does not divide $\zeta$, it does not divide any 
$\zeta_i$. Write
$$
c(X,Y) = c_1 \left( X^{1/m}, Y^{1/m} \right) ... c_s \left( X^{1/m}, Y^{1/m}
\right),
$$
the factorization of $c(X,Y)$ in $S$. There must be some $c_i$ not dividing
$\zeta$, so write $\zeta_1,...,\zeta_j$ the conjugates of $\zeta$ not divided
by $c_i$. Then
$$
a_{n-j} (X,Y) = \sum_{l_1,...,l_j \in \{1,...,n\}} \zeta_{l_1}...
\zeta_{l_j},
$$
so $c_i$ divides all terms of the right hand side except $\zeta_1...\zeta_j$,
which is imposible as $c(X,Y) | a_{n-j}$.

The only remaining posibility is that $c_i$ does not divide any conjugate of
$\zeta$. This is clearly absurd again, as $c_i$ must divide $a_0 = \zeta_1...
\zeta_n$.

\subsection{The non--transversal case}

Write the minimal polynomial of $\zeta \in S$ as
$$
F(Z) = Z^n + ... + a_0(X,Y),
$$
with $\nu (F) = l < n$. We will look up first for equimultiple curves $P$, not 
lying in $Z=0$. If the initial form of $\zeta$ contains both $X^{1/m}$ and $Y^{1/m}$
it comes out that both $X$ and $Y$ appear in $\ol{F}$. Then, as $D(F) \in P$, for 
any derivation $D$ of order $l-1$, it is plain that $P = (X,Y)$, which is impossible
since $F$ is a monic polynomial in $Z$.

Hence $\zeta$ can be assumed to be
$$
\zeta = X^{a/m} + \zeta',
$$
with $a<m$, $\nu(\zeta') > \nu (\zeta)$. But, $\zeta'$ contains a monomial of the form
$X^{b/m}Y^{c/m}$, with $b<a$, and we take one of these monomials with minimal $b$,
it follows (for the same reason as above) that $P=(X,Y)$.

So, if there are equimultiple curves not lying in $Z=0$, the series has the form 
$\zeta = X^\lambda u(X^{1/m},Y^{1/m})$, and, therefore, $(X,Z)$ is an equimutiple
curve as well. This situation features a curious phenomenon, which is shown in
the following lemma, whose proof follows the same lines as Lipman's
normalization for quasi--ordinary surfaces (\cite{Lipman}).

\noindent {\bf Lemma.--}
Let $\zeta \in S$, $\nu (\zeta)<1$ and $(X,Z)$ lying in the equimultiple locus of
the surface $\cS$, defined by $\zeta$. Then there is a transversal Puiseux 
power series $\eta \in S$ defining $\cS$.

\noindent {\bf Proof:}
Let us write, as above, $\zeta = X^{\lambda} u(X^{1/m}, Y^{1/m})$. Let us
call $F$ the minimal polynomial of $\zeta$; $\mu/m$ the order of 
$\zeta$. $F$ is regular with respect to $Z$ and $X$ (with order
$\lambda m^2/d$). So, let us call $\widehat{F}$ the Weierstrass polynomial 
associated to $F$ with respect to $X$.

Let now $T$ be a new variable. As the units of $R$ have all of their $\mu$--th
roots inside $R$, the power series
$$
X^\mu u(X,Y) - T^\mu \in \CC[[X,Y,T]]
$$
has a factor $Xu'(X,Y)-T$, where $u'$ is an $\mu$--th root of $u$. 
By the Weierstrass Preparation Theorem, there exists a unit $u''(X,Y,T)$
such that
$$
u''(X,Y,T)(Xu'(X,Y)-T) = X - R(Y,T) \quad \quad \quad \quad (\star)
$$

Doing $X=0$ we can easily find $R(Y,T) = Tu''(0,Y,T)$.
Define then $G(Y,T) = u''(0,Y,T)$, which is a unit in 
$\CC[[Y,T]]$. Writing $X = R(Y,T)$ in the above formula we get
$$
0 = u''(R(Y,T),Y,T) \left[ R(Y,T)u'(R(Y,T),Y)-T \right].
$$
Hence, as $R(Y,T) = TG(Y,T)$,
$$
TG(Y,T)u'(TG(Y,T),Y) - T = 0.
$$

This means $X = TG(Y,T)$ annihilates a factor of $(\star)$ and therefore
$$
\left(TG(Y,T)\right)^\mu u(TG(Y,T),Y) - T^\mu = 0.
$$

If we write now $X^{1/m}$ instead of $X$, $Y^{1/m}$ instead of $Y$ and
$Z^{1/\mu}$ instead of $T$, we can assure that there is a unit
$G(Y^{1/m},Z^{1/\mu})$ such that, if we write $X^{1/m} =
Z^{1/\mu}G(Y^{1/m},Z^{1/\mu})$, we get
$$
\zeta \left( Z^{1/\mu}G \left( Y^{1/m},Z^{1/\mu} \right), Y^{1/m}
\right) = Z.
$$

As we have
$$
F\left( X,Y,\zeta \left( X^{1/m},Y^{1/m} \right) \right)=0,
$$
doing $X^{1/m} = Z^{1/\mu}G(Z^{1/\mu},Y^{1/m})$ we obtain
$$
F \left( Z^{m/\mu} G^m \left( Z^{1/\mu},Y^{1/m} \right),Y,Z \right) = 0.
$$

Hence $\eta = Z^{m/\mu} G^m \left( Z^{1/\mu},Y^{1/m} \right)$ is a transversal
Puiseux power series which is a root of $\widehat{F}$.

\noindent {\bf Remark.--}
As for other kind of primes in the equimultiple locus of $\cS$, note that the 
arguments for the transversal case do not use transversality at all. That
means $(Z,c(X,Y))$ lies in the equimultiple locus if and only if $c(X,Y) | \zeta$,
which is imposible when $\nu (\zeta) < 1$.

This finishes the considerations for the non--transversal case.

We summarize the results on this section.

\noindent {\bf Proposition.--}
Let $\cS$ be a Puiseux surface defined by $\zeta \in S$ and $F(Z)$ the minimal
polynomial of $\zeta$. Then

\begin{enumerate}
\item[(a)] If $\nu(\zeta)\geq1$, a prime ideal $(Z,c(X,Y))$ lies in the equimultiple
locus if and only if $c(X,Y) | \zeta$.
\item[(b)] If $\nu(\zeta)<1$ either there are no primes of height $1$ in the equimultiple
locus or there exists a change of variables which takes the surface into the 
transversal case.
\end{enumerate}

\noindent {\bf Remark.--}
Again, it is not posible to express the equimultiple locus from
a set of exponents of $\zeta$, as it is shown in the following example:
$$
\zeta_1 = X^{3/2} + Y^{1/2}X - YX^{1/2} - Y^{3/2}, \; \;
\zeta_2 = X^{3/2} + 3Y^{1/2}X + 3YX^{1/2} + Y^{3/2}.
$$

The first Puiseux power series defines a surface with equation
$$
F_1 = Z^4 - 2(X-Y)^2(X+Y)  Z^2 + (X-Y)^6
$$
where the prime $P=(Z,X-Y)$ lies in the equimultiple locus, while an equation for
the second surface is
$$
F_2 = Z^4 - 2 \left( X^3 + 15 X^2Y + 15XY^2 + Y^3 \right) Z^2 + (X-Y)^6,
$$
hence the equimultiple locus of this surface is simply the closed point.

\section{Final comments}

Puiseux singularities are still far from being well understood. In particular,
two goals seem highly interesting to attain:

\begin{enumerate}
\item[(a)] A simple (as simple as possible, at least) criterion for deciding
whenever a given surface is Puiseux. In the quasi--ordinary case, this is the
celebrated Jung--Abhyankar theorem (\cite{Jung}, \cite{Abhy}) but, for the
more general case, we cannot know whether a given polynomial has Puiseux roots
unless we actually compute them. This is now posible, thanks to the work of
McDonald (\cite{McD}) and Aroca--Cano (\cite{Fuen}); but, while useful for
working with examples, this method has not been adopted so far for its use on
abstract argumentations.
\item[(b)] A numerical (say arithmetical, say combinatorial) control of the
behaviour of the singularity in the resolution process. In the quasi--ordinary
case, this is the main result in Lipman's thesis (\cite{Lipman}). For the
general case, only partial results are known (see \cite{PhD}), which we 
summarize now:

(i) It is possible to describe a resolution process for Puiseux surfaces which,
in each step, preserves the Puiseux character. This resolution process is,
basicly, the Levi resolution (\cite{BL}) which consists in blowing--up maximal
subvarieties of the equimultiple locus; up to some previous monoidal 
transformations which take our surface into another one whose roots are
$\nu$--quasi--ordinary (this is a mild generalization of the quasi--ordinary 
concept, see, for instance, \cite{Luengo}).

(ii) For some cases (transversal and $\nu$--quasi--ordinary) it is possible
to prove that the blowing--up of a Puiseux surface is still a Puiseux surface,
and also we are able to know a set of privileged exponents of the
transformed surface from a set of the original one. Still, we do not know
what happens in the general case.
\end{enumerate}

\noindent {\bf Example.--} Let us end with an example showing the difficulties of the 
behaviour under blowing--up of a Puisuex surface. The computations below have
been posible thanks to the software developed by F. Aroca and J. Cano (\cite{Fuen}).

Let us consider 
$$
\zeta = X^{2/5} + Y^{2/5},
$$
whose minimal polynomial is
$$
F(Z) = Z^{25} + a_{20}(X,Y) Z^{20} + a_{15}(X,Y) Z^{15} + a_{10}(X,Y) Z^{10} +
a_{5}(X,Y) Z^{5} + a_{0}(X,Y),
$$
where
$$
\begin{array}{rcl}
a_{20}(X,Y) & = & -5X^2-5Y^2 \\
a_{15}(X,Y) & = & 10X^4+10Y^4-605X^2Y^2 \\
a_{10}(X,Y) & = & -10X^6-10Y^6-1905X^4Y^2-1905X^2Y^4 \\
a_5 (X,Y) & = & 5X^8+5Y^8-605X^6Y^2-605X^2Y^6+1905X^4Y^4 \\
a_0 (X,Y) & = & -X^{10}-Y^{10}-5X^8Y^2-5X^2Y^8-10X^4Y^6-10X^6Y^4
\end{array}
$$

A quadratic transform in the point $(0:0:1)$ of the tangent cone gives the
equation
$$
Z^{15} + a_{20}(X,Y) Z^{12} + a_{15} (X,Y) Z^9 + a_{10}
(X,Y) Z^6 + a_5 (X,Y) Z^3 + a_0 (X,Y).
$$

The roots of this equation (which is a Weierstrass equation w.r.t. $X$, $Y$
and $Z$) as a polynomial in $Z$ are
$$
\eta = Y^{2/3} + \frac{5}{3} X^{2/5}Y^{4/15} + \frac{5}{9}
\frac{X^{4/5}}{Y^{2/15}} - \frac{5}{81} \frac{X^{6/5}}{Y^{8/15}};
$$
and all of its conjugates. Analogously, its roots as a polynomial in $X$ are
$$
\varsigma = Z^{3/2} - \frac{5}{2} Z^{9/10}Y^{2/5} + \frac{15}{8}
Z^{3/10}Y^{4/5} - \frac{5}{16} \frac{Y^{6/5}}{Z^{3/10}},
$$
and all of its conjugates. As the polynomial is symmetric in $X$ and $Y$, we
cannot find an easy parametrization which gives us an equation with Puiseux
roots. The existence of such an equation is a problem which relates with the
two open questions set above.

\end{document}